# Best proximity points and weak convergence theorems for hybrid mappings in Banach spaces

M. De la Sen

Instituto de Investigacion y Desarrollo de Procesos. Universidad del Pais Vasco

Campus of Leioa (Bizkaia) – Aptdo. 644- Bilbao, 48080- Bilbao. SPAIN

**Abstract**: The paper gives some results on best proximity and fixed point for a class of generalized hybrid cyclic self-mappings in Banach spaces.

## 1. Introduction

The following objects are considered through the manuscript:

1) The Hilbert space $H$ on the field $X$ (in particular, $\mathbf{R}$ or $\mathbf{C}$) endowed with the inner product $\langle x, y \rangle$ which maps $H \times H$ to $X$; $\forall x, y \in H$ which maps $H \times H$ to $X$; where $(X, \|\ \|)$ is a Banach space when endowed with a norm $\|\ \|$ induced by the inner product and defined by $\|x\| = \langle x, x \rangle^{1/2}$; $\forall x \in H$. It is well-known that all Hilbert spaces are uniformly convex Banach spaces and that Banach spaces are always reflexive.

2) The $p(\geq 2)$-cyclic self-mapping $T: A \to A$ with $A := \bigcup_{i \in \bar{p}} A_i$ subject to $A_{p+1} \equiv A_p$ where $A_i (\neq \varnothing) \subset H$ are $p$ subsets of $H$; $\forall i \in \bar{p} = \{1, 2, \ldots, p\}$, that is, a self-mapping satisfying $T(A_i) \subseteq A_{i+1}$; $\forall i \in \bar{p}$

3) The function $f: D (\equiv dom\ f) \subset X \to (-\infty, \infty]$ which is a proper convex function which is Gâteaux differentiable in the topological interior of $D$, $int\ D$, that is, $D := \{x \in X : f(x) < \infty\} \neq \varnothing$ and convex since $f$ is proper with

$$f(\alpha x + (1-\alpha)y) \leq \alpha f(x) + (1-\alpha)f(y); \quad \forall x, y \in D, \ \forall \alpha \in [0, 1] \tag{1}$$

since $f$ is convex, and for each $x \in D$, there is $x^* = f(x) \in X^*$ (the topological dual of $X$) such that

$$\exists \lim_{t \to 0} \frac{f(x+ty) - f(x)}{t} = \langle y, f'(x) \rangle; \quad \forall y \in D \tag{2}$$

since $f$ is Gâteaux differentiable in $int\ D$ where $f'(x)$ denotes the Gâteaux derivative of $f$ at $x$ if $x \in int\ D$. On the other hand, $f$ is said to be strictly convex if

$$f(\alpha x + (1-\alpha)y) < \alpha f(x) + (1-\alpha)f(y);\ \forall x, y(\neq x) \in D,\ \forall \alpha \in (0,1) \tag{3}$$

4) The Bregman distance (or Bregman divergence) $D_f$ associated with the proper convex function $f$

$D_f : D \times D \to (-\infty, \infty]$, where $\mathbf{R}_{0+} := \{z \in \mathbf{R} : z \geq 0\} = \mathbf{R}_+ \cup \{0\}$, is defined by:

$$D_f(y,x) = f(y) - f(x) - \langle y - x, f'(x) \rangle;\ \forall x, y \in D \tag{4}$$

provided that it is Gâteaux differentiable everywhere in $int\ D$. If $f$ is not Gâteaux differentiable at $x \in int\ D$ then (4) is replaced by

$$D_f(y,x) = f(y) - f(x) + f^0(x, x-y) \tag{5}$$

where $f^0(x, x-y) := \lim_{t \to 0^+} \dfrac{f(x + t(x-y)) - f(x)}{t}$ and $D_f(y,x)$ is finite if and only if $x \in D^0 \subset D$, the algebraic interior of $D$ defined by:

$$D^0 := \{x \in D : \exists z \in (x,y), [x,z] \subseteq D;\ \forall y \in X \setminus \{x\}\} \tag{6}$$

The topological interior of $D$ is $int(D) := \{x \in D : x \in fr(D)\} \subset D^0$. where $fr(D)$ is the boundary of $D$. It is well-known that the Bregman distance does not satisfy either the symmetry property or the triangle inequality which are required for standard distances while they are always nonnegative because of the convexity of the function $f : D(\equiv dom\ f) \subset X \to (-\infty, \infty]$. The Bregman distance between sets $B, C \subset H \subset X$ is defined as $D_f(B,C) := \inf\limits_{x \in B,\ y \in C} D_f(x,y)$. If $A_i \in A \subset H$ for $i \in \overline{p}$ then

$D_{fi} := D_f(A_i, A_{i+1}) = \inf\limits_{x \in A_i,\ y \in A_{i+1}} D_f(x,y)$. Through the paper, sequences $\{x_n\}_{n \in N_0} \equiv \{T^n x\}_{n \in N_0}$ with $N_0 = N \cup \{0\}$ are simply denoted by $\{T^n x\}$ for the sake of notation simplicity.

**2. Some fixed point theorems for cyclic hybrid self-mappings on the union of intersecting subsets**

The Bregman distance is not properly a distance, since it does not satisfy symmetry and the triangle inequality, but it is always nonnegative and leads to the following interesting result towards its use in applications of fixed point theory:

**Lemma 2.1**. If $f : D \times D \to (-\infty, \infty]$ is a proper strictly convex function being Gâteaux differentiable in *int D* then

$$D_f(y, x) + D_f(x, y) = -\langle y - x, f'(x) - f'(y)\rangle \geq 0 ; \quad \forall x, y \in int\ D \tag{7}$$

$$D_f(x, x) = 0 ; \quad \forall x \in int\ D \tag{8}$$

$$D_f(y, x) > 0 ; \forall x, y(\neq x) \in int\ D \tag{9}$$

$$D_f(y, x) - D_f(x, y) = 2(f(y) - f(x)) - \langle y - x, f'(x) + f'(y)\rangle ; \forall x, y \in int\ D \tag{10}$$

*Proof*: By using (4) for $D_f(y, x)$ and defining $D_f(x, y) = f(x) - f(y) - \langle x - y, f'(y)\rangle$ ; $\forall x, y \in int\ D$ by interchanging $x$ and $y$ in the definition of $D_f(y, x)$ in (4)

$$D_f(y, x) + D_f(x, y) = \langle x, f'(x)\rangle + \langle y, f'(y)\rangle - \langle x, f'(y)\rangle > -\langle y, f'(x)\rangle$$

$$= \langle x, f'(x) - f'(y)\rangle + \langle y, f'(y) - f'(x)\rangle \tag{11}$$

what leads to the equality in (7). The whole (7) follows from (11), the fact that $\langle x - y, f'(x) - f'(y)\rangle > 0$ ; $\forall x, y(\neq x) \in int\ D$, [5], if $f : D \times D \to (-\infty, \infty]$ is proper strictly convex, and the fact $D_f(x, y) \geq 0$ ; $\forall x, y \in int\ D$.

Eq. (8) follows from (7) for $x = y$ leading to $2D_f(x, x) = 0$. To prove (9), take $x, y(\neq x) \in int\ D$ implying that $\langle x - y, f'(x) - f'(y)\rangle > 0$, [5], and proceed by contradiction using (4) by assuming that $D_f(y, x) = 0$ for such $x, y(\neq x) \in int\ D$ so that

$$0 = D_f(y, x) = f(y) - f(x) - \langle y - x, f'(x)\rangle$$

$$= f(y) - f(x) + \langle x - y, f'(x) - f'(y)\rangle + \langle x - y, f'(y)\rangle$$

$$> f(y) - f(x) - \langle y - x, f'(y)\rangle = D_f(y, x)$$

which contradicts $D_f(y, x) = 0$. Then, $D_f(y, x) > 0$ and hence (9).

$$D_f(y, x) - D_f(x, y) = f(y) - f(x) + \langle x - y, f'(x)\rangle - f(x) + f(y) + \langle x - y, f'(y)\rangle$$

$$= 2(f(y) - f(x)) + \langle x - y, f'(y) + f'(x)\rangle \tag{12}$$

; $\forall x, y \in int\ D$ and hence (10)-(11) via (7) and (9). □

The following definition is then used:

**Definition 2.2**. If $D \cap A_i \neq \emptyset$; $\forall i \in \bar{p}$ and $f : D(\equiv dom\, f) \subset X \to (-\infty, \infty]$ is a proper convex function which is Gâteaux differentiable in $int\, D$, then the $p$-cyclic self-mapping $T : \bigcup_{i \in \bar{p}} A_i \to \bigcup_{i \in \bar{p}} A_i$, where $A := \bigcup_{i \in \bar{p}} A_i \subseteq int\, D \subset H$ and $A_i \neq \emptyset$; $\forall i \in \bar{p}$, is said to be a *generalized* point-dependent $(K, \lambda)$-hybrid $p(\geq 2)$-cyclic self-mapping relative to $D_f$ if

$$D_f(Tx, Ty) \leq K_i(y) D_f(x, y) + \lambda(y)\langle x - Tx, f'(y) - f'(Ty)\rangle \; ; \; \forall x \in A_i, \forall y \in A_{i+1}, \forall i \in \bar{p} \quad (13)$$

for some given functions $\lambda : \bigcup_{i \in \bar{p}} A_i \to \mathbf{R}$ and $K_i : A_{i+1} \to (0, a_i]$ with $a_i \in \mathbf{R}_+$, $\forall i \in \bar{p}$, where $K : \bigcup_{i \in \bar{p}} A_{i+1} \to (0, 1]$ defined by $K(y) = \prod_{j=i}^{i+p-1}\left[K_j(T^{j-i}y)\right]$ for any $y \in A_{i+1}$; $\forall i \in \bar{p}$.

If, furthermore, $K : \bigcup_{i \in \bar{p}} A_i \to (0, 1)$; $\forall i \in \bar{p}$ then $T : \bigcup_{i \in \bar{p}} A_i \to \bigcup_{i \in \bar{p}} A_i$ is said to be a generalized point-dependent $(K, \lambda)$-hybrid $p(\geq 2)$-cyclic self-mapping relative to $D_f$. □

If $p = 1$, it is possible to characterize $T : A_1 \to A_1$ as a trivial 1-cyclic self-mapping with $A_1 = A_2$ which does not need to be specifically referred to as 1-cyclic.

The following concepts are useful:

$f : D(\equiv dom\, f) \subset X \to (-\infty, \infty]$ is said to be totally convex if the modulus of total convexity $v_f : D^0 \times [0, \infty) \to [0, \infty]$, that is, $v_f(x, t) = inf\{D_f(x, y) : y \in D, \|y - x\| = t\}$ is positive for $t > 0$.

$f : D(\equiv dom\, f) \subset X \to (-\infty, \infty]$ is said to be uniformly convex if the modulus of total convexity $\delta_f : [0, \infty) \to [0, \infty]$, that is, $\delta_f(t) = inf\left\{f(x) + f(y) - 2f\left(\frac{x+y}{2}\right) : x, y \in D, \|y - x\| = t\right\}$ is positive for $t > 0$. It holds that $v_f(x, t) \geq \delta_f(t)$; $\forall x \in D$, [5]. The following result holds:

**Theorem 2.3**. Assume that:

1) $f : D(\equiv dom\, f) \subset X \to (-\infty, \infty]$ is a lower-semicontinuous proper strictly totally convex function which is Gâteaux differentiable in $int\, D$.

2) $A_i(\neq \emptyset) \subseteq int\, D \subset H$; $\forall i \in \bar{p}$ are bounded, closed and convex subsets of $H$ which intersect and $T : \bigcup_{i \in \bar{p}} A_i \to \bigcup_{i \in \bar{p}} A_i$ is a generalized point-dependent $(K, \lambda)$-hybrid $p(\geq 2)$-cyclic self-mapping relative to $D_f$ for some given functions $\lambda : \bigcup_{i \in \bar{p}} A_i \to \Lambda \subset \mathbf{R}$ and $K : \bigcup_{i \in \bar{p}} A_i \to (0, 1)$, defined by $K(y) = \prod_{j=i}^{i+p-1}\left[K_j\left(T^{j-i}y\right)\right]$ for any $y \in A_{i+1}$; $\forall i \in \bar{p}$ and some functions $K_i : A_{i+1} \to (0, a_i]$; $\forall i \in \bar{p}$, with $\Lambda$ being bounded.

3) There is a convergent sequence $\{T^n x\}$ to some $z \in \bigcap_{i \in \bar{p}} A_i$ for some $x \in \bigcup_{i \in \bar{p}} A_i$.

Then, $z = Tz$ is the unique fixed point of $T : \bigcup_{i \in \bar{p}} A_i \to \bigcup_{i \in \bar{p}} A_i$ to which all sequences $\{T^n x\}$ converge for any $x \in \bigcup_{i \in \bar{p}} A_i$; $\forall i \in \bar{p}$.

**Proof**: The recursive use of (13) yields:

$$D_f\left(T^2 x, T^2 y\right) \leq K_{i+1}(Ty) D_f(Tx, Ty) + \lambda(Ty)\langle Tx - T^2 x, f'(Ty) - f'(T^2 y)\rangle$$

$$\leq K_{i+1}(Ty)\left[K_i(y) D_f(x, y) + \lambda(y)\langle x - Tx, f'(y) - f'(Ty)\rangle\right] + \lambda(Ty)\langle Tx - T^2 x, f'(Ty) - f'(T^2 y)\rangle$$

; $\forall x \in A_i, \forall y \in A_{i+1}, \forall i \in \bar{p}$ (14)

..................................................................................

$$D_f\left(T^p x, T^p y\right) \leq K_{i+p-1}\left(T^{p-1} y\right) D_f\left(T^{p-1} x, T^{p-1} y\right) + \lambda\left(T^{p-1} y\right)\langle T^{p-1} x - T^p x, f'\left(T^{p-1} y\right) - f'\left(T^p y\right)\rangle$$

$$\leq \left[\prod_{j=1}^{p} K_{p-j+1}\left(T^{p-j+1-i} y\right)\right] D_f(x, y)$$

$$+ \sum_{k=1}^{p}\left(\prod_{j=k+1}^{p}\left[K_{p-j+i}\left(T^{p-j+1}\right)y\right]\right)\lambda\left(T^{k-1} y\right)\langle T^{k-1} x - T^k x, f'\left(T^{k-1} y\right) - f'\left(T^k y\right)\rangle \quad (15)$$

; $\forall x \in A_i, \forall y \in A_{i+1}, \forall i \in \bar{p}$ with $T^p x \in A_{i+p}$, $T^p y \in A_{i+1+p}$ with $A_{i+p} = A_i$, $K_{i+p} = K_i$; $\forall i \in \bar{p}$.

where $T^0$ is the identity mapping on $\bigcup_{i \in \bar{p}} A_i$. Now, define $\hat{K}(y) := \left[\prod_{j=1}^{p} K_{p-j+1}\left(T^{p-j+1-i} y\right)\right]$ so that

one gets

$$D_f\left(T^{np} x, T^{np} y\right) \leq \hat{K}^n(y) D_f(x, y) + \sum_{k=1}^{np}\left(\prod_{j=k+1}^{np}\left[K_{np-j+i}\left(T^{np-j+1}\right)y\right]\right)\lambda\left(T^{k-1} y\right)\langle T^{k-1} x - T^k x, f'\left(T^{k-1} y\right) - f'\left(T^k y\right)\rangle$$

$$\leq \hat{K}^n(y) D_f(x, y) + \sum_{k=1}^{(n-1)p}\left(\prod_{j=k+1}^{np}\left[K_{np-j+i}\left(T^{np-j+1}\right)y\right]\right)\lambda\left(T^{k-1} y\right)\langle T^{k-1} x - T^k x, f'\left(T^{k-1} y\right) - f'\left(T^k y\right)\rangle$$

$$+ \sum_{k=(n-1)p}^{np}\left(\prod_{j=k+1}^{np}\left[K_{np-j+i}\left(T^{np-j+1}\right)y\right]\right)\lambda\left(T^{k-1} y\right)\langle T^{k-1} x - T^k x, f'\left(T^{k-1} y\right) - f'\left(T^k y\right)\rangle$$

$$\leq \hat{K}^n(y)D_f(x,y) + \left(\frac{1-\hat{K}^{(n-1)p+1}(\hat{y})}{1-\hat{K}(y)} + M_{np}\right)\left(\max_{1\leq j\leq np}\left[\lambda(T^{j-1}y)\langle T^{j-1}x - T^j x, f'(T^{j-1}y) - f'(T^j y)\rangle\right]\right)$$

(16)

since $\hat{K}(y) = \prod_{j=i}^{i+p-1}\left[K_j(T^{j-i}y)\right] < 1$; $\forall y \in A_{i+1}$ since $T: \bigcup_{i\in\bar{p}} A_i \to \bigcup_{i\in\bar{p}} A_i$ is a generalized point-dependent $(K,\lambda)$-hybrid $p(\geq 2)$-cyclic self-mapping relative to $D_f$, implies that $\hat{K}(y) < 1$; $\forall y \in \bigcup_{i\in\bar{p}} A_i$, where

$$M_{np} \geq \sum_{k=(n-1)p}^{np}\left(\prod_{j=k+1}^{np}\left[K_{np-j+i}(T^{np-j+1})y\right]\right)\lambda(T^{k-1}y)\langle T^{k-1}x - T^k x, f'(T^{k-1}y) - f'(T^k y)\rangle$$

(17)

and

$$D_f(T^{mnp}x, T^{mnp}y) \leq$$
$$\hat{K}^m(y)D_f(T^{np}x, T^{np}y) + \left(\frac{1-\hat{K}^{(m-1)p+1}(\hat{y})}{1-\hat{K}(y)} + M_{nmp}\right)\left(\max_{np+1\leq j\leq nmp}\left[\lambda(T^{j-1}y)\langle T^{j-1}x - T^j x, f'(T^{j-1}y) - f'(T^j y)\rangle\right]\right)$$

(18)

so that

$$0 \leq \limsup_{n,m\to\infty} D_f(T^{mnp}x, T^{mnp}y) \leq \left(\frac{1}{1-\hat{K}(y)} + \limsup_{n,m\to\infty} M_{nmp}\right)$$
$$\times \limsup_{n,m\to\infty}\left(\max_{np+1\leq j\leq nmp}\left[\lambda(T^{j-1}y)\langle T^{j-1}x - T^j x, f'(T^{j-1}y) - f'(T^j y)\rangle\right]\right) = 0 \quad (19)$$

since, $\lambda: \bigcup_{i\in\bar{p}} A_i \to \Lambda \subset \mathbf{R}$ is bounded, $f: D \subset X \to (-\infty, \infty]$ is lower-semicontinuous then with all subgradients in any bounded subsets of $int\, D$ being bounded, and $\{T^j x\}$ and $\{T^{j-1}x - T^j x\}$; $\forall x \in \bigcup_{i\in\bar{p}} A_i, \forall i \in \bar{p}$ converge so that they are Cauchy sequences being then bounded; $\forall x \in \bigcup_{i\in\bar{p}} A_i, \forall i \in \bar{p}$ where $z \in \bigcap_{i\in\bar{p}} A_i$, since $\bigcap_{i\in\bar{p}} A_i$ is nonempty and closed, is some fixed point of $T: \bigcup_{i\in\bar{p}} A_i \to \bigcup_{i\in\bar{p}} A_i$. As a result, $\exists \lim_{n\to\infty} D_f(T^{np}x, T^{np}y) = \lim_{n\to\infty} D_f(T^n x, T^n y) = 0$; $\forall x \in A_i, \forall y \in A_{i+1}, \forall i \in \bar{p}$. From a basic property of Bregman distance, $T^n y \to T^n x (\to z)$ as $n\to\infty$; $\forall x \in A_i, \forall y \in A_{i+1}, \forall i \in \bar{p}$, if $f: D(\equiv dom\, f) \subset X \to (-\infty, \infty]$ is sequentially consistent. But, since $\bigcup_{i\in\bar{p}} A_i$ is closed, $f: D|\bigcup_{i\in\bar{p}} A_i \to \bigcup_{i\in\bar{p}} A_i$ is sequentially consistent if and only if it is totally convex, [6]. Thus, $\{T^n y\}$ converges also to $z$ for any $x \in A_i$ and $y \in A_{i+1}$; $\forall i \in \bar{p}$ so that $z = Tz$ is a fixed point of $T: \bigcup_{i \in \bar{p}} A_i \to \bigcup_{i \in \bar{p}} A_i$. Assume not and proceed by contradiction leading to $D_f(T^n x, T^n z) \to D_f(z, T^n z)$ converges to zero as $n \to \infty$ from a basic property of Bregman distance. Thus, $\left[ D_f(z, T^n z) - f(z) + f(T^n z) \right] \to 0$ as $n \to \infty$ since $\langle z - T^n z, f'(T^n z) \rangle \to 0$ as $n \to \infty$. As a result, $f(T^n z) \to f(z)$, $T^n z \to z = Tz$ as $n \to \infty$ from the continuity of $f: D(\equiv dom\, f) \subset X \to (-\infty, \infty]$ and $z$ is a fixed point of $T: \bigcup_{i \in \bar{p}} A_i \to \bigcup_{i \in \bar{p}} A_i$. Now, take any $y_1 \in A_j$ so that $y = T^{i+1-j} y_1 \in A_{i+1}$ then

$D_f(T^n x, T^{n+i+1-j} y) \to D_f(z, T^{i+1-j} z) = D_f(z, z) = 0$ since $z$ is a fixed point of $T: \bigcup_{i \in \bar{p}} A_i \to \bigcup_{i \in \bar{p}} A_i$ and $f: D(\equiv dom\, f) \subset X \to (-\infty, \infty]$ is a proper strictly totally convex function. As a result, $\{T^n y\}$ converges to $z$; $\forall y \in A_{i+1}$.

It is now proven that $z \in \bigcap_{i \in \bar{p}} A_i$ is the unique fixed point of $T: \bigcup_{i \in \bar{p}} A_i \to \bigcup_{i \in \bar{p}} A_i$. Assume not so that there is $z_1 (\neq z) = T z_1 \in \bigcap_{i \in \bar{p}} A_i$. Then, $D_f(T^n x, T^n z_1) \to 0$ as $n \to \infty$ so that $\{T^n x\} \to z_1 = T^n z_1$, since $D_f(x, y) > 0$ if $x, y (\neq x) \in \bigcap_{i \in \bar{p}} A_i = int\left(\bigcap_{i \in \bar{p}} A_i\right)$ since $f: D(\equiv dom\, f) \subset X \to (-\infty, \infty]$ is proper and totally strictly convex and, since, $\{T^n x\} \to z$, and $z \in \bigcap_{i \in \bar{p}} A_i$. Since $\bigcap_{i \in \bar{p}} A_i$ is closed and convex, it turns out that $z$ is the unique fixed point of $T: \bigcup_{i \in \bar{p}} A_i \to \bigcup_{i \in \bar{p}} A_i$.

Note that the result also holds for any $\forall y_1 \in \bigcup_{i \in \bar{p}} A_i$ since $y \in \bigcup_{j(\neq i) \in \bar{p}} A_j$ maps to $y_1 = T^{k_i} y \in A_{i+1}$ for some nonnegative integer $k_i \leq p - 1$ through the self-mapping $T: \bigcup_{i \in \bar{p}} A_i \to \bigcup_{i \in \bar{p}} A_i$ so that $D_f(T^n x, T^n y_1) = D_f(T^n x, T^{n+k_i} y) \to D_f(z, T^{k_i} z) = D_f(z, z) = 0$ as $n \to \infty$ since $z \in \bigcap_{i \in \bar{p}} A_i$ is the unique fixed point of $T: \bigcup_{i \in \bar{p}} A_i \to \bigcup_{i \in \bar{p}} A_i$ and $\{T^n y\}$ converges to $z$ for any $y \in \bigcup_{i \in \bar{p}} A_i$. □

The subsequent result extends directly Theorem 2.3 to the $p$-composite self-mappings $T_i^p: \bigcup_{j \in \bar{p}} A_j | A_i \to A_i$; $\forall i \in \bar{p}$, defined as $T^p x = T^j (T^{p-j} x)$; $\forall x \in \bigcup_{i \in \bar{p}} A_i$ subject to $i = p - j - k$;

$\forall i, j \in \bar{p}$. The subsets $A_i \subset X$, $i \in \bar{p}$ are not required to intersect since the restricted composite mappings as defined above are self-mappings on nonempty, closed and convex sets.

**Corollary 2.4.** Assume that:

1) $f_i : D(\equiv dom\, f) \subset X \to (-\infty, \infty]$ is a proper strictly totally convex function which is lower-semicontinuous and Gâteaux differentiable in $int\, D$.

2) $A_i(\neq \varnothing) \subseteq int\, D \subset H$ is bounded and closed; $\forall i \in \bar{p}$, $T : \bigcup_{i \in \bar{p}} A_i \to \bigcup_{i \in \bar{p}} A_i$ is a $p$-cyclic self-mapping so that $T_i^p : \bigcup_{j \in \bar{p}} A_j | A_i \to A_i$ for some $i \in \bar{p}$ is a generalized point-dependent $(K, \lambda_i)$-hybrid $p(\geq 2)$-cyclic self-mapping relative to $D_f$ for some given functions $\lambda_i : \bigcup_{i \in \bar{p}} A_i \to \Lambda \subset \mathbf{R}$ and $K : \bigcup_{i \in \bar{p}} A_i \to (0,1)$ for some $i \in \bar{p}$ defined by $K(y) = \prod_{j=i}^{i+p-1}\left[K_j\left(T^{j-i}y\right)\right]$ for any $y \in A_{i+1}$; $\forall i \in \bar{p}$ and $K_i : A_{i+1} \to (0, a_i]$ for some $a_i \in \mathbf{R}_+$, $\forall i \in \bar{p}$ where $A := \bigcup_{i \in \bar{p}} A_i \subset H$, $\Lambda$ being bounded and $A_i$ being, furthermore, convex for the given $i \in \bar{p}$.

3) There is a convergent sequence $\{T_i^{np} x\}$ to some $z_i \in A_i$ for some $x \in A_i$ and $i \in \bar{p}$.

Then, $z_i = Tz_i$ is a unique fixed point of $T_i^p : \bigcup_{j \in \bar{p}} A_j | A_i \to A_i$ to which all sequences $\{T_i^{np} x\}$ converge for any $x \in A_i$ for $i \in \bar{p}$.

Also, if Conditions 1-3 are satisfied with all the subsets $A_i$; $\forall i \in \bar{p}$ being nonempty, closed and convex for some proper strictly convex function $f \equiv f_i : D(\equiv dom\, f) \subset X \to (-\infty, \infty]$ which is Gâteaux differentiable in $int\, D$, then $z_i = Tz_i$; $\forall i \in \bar{p}$ is a unique fixed point of $T_i^p : \bigcup_{j \in \bar{p}} A_j | A_i \to A_i$; $\forall i \in \bar{p}$ to which all sequences $\{T_i^{np} x\}$ converge for any $x \in A_i$; $\forall i \in \bar{p}$. The $p$ unique fixed points of each generalized point-dependent $(K, \lambda_i)$-hybrid 1-cyclic composite self-mappings $T_i^p : \bigcup_{j \in \bar{p}} A_j | A_i \to A_i$; $\forall i \in \bar{p}$ fulfil the relations $z_{p-i} = T^j z_k$ for $i = p - j - k$; $\forall i \in \overline{p-1}$, $\forall j \in \bar{p}$. □

**Outline of proof**: Note that $D \cap \bigcap_{i \in \bar{p}}(A_i) \neq \varnothing$. Eq. 13 is now extended to $T_i^p : \bigcup_{j \in \bar{p}} A_j | A_i \to A_i$ for the given $i \in \bar{p}$ leading to

$$D_f\left(T_i^p x, T_i^p y\right) \leq K_i(y) D_f(x, y) + \lambda_i(y)\langle x - T_i^p x, f'(y) - f'(T_i^p y)\rangle \ ; \ \forall x, y \in A_i, i \in \bar{p} \qquad (20)$$

since $T_i^p$ is a trivial $1$-cyclic self-mapping on $A_i$ for $i \in \bar{p}$. The above relation leads recursively to:

$$D_f\left(T_i^{np} x, T_i^{np} y\right) \leq \hat{K}^n(y) D_f(x, y)$$

$$+ \left(\frac{1}{1-\hat{K}(y)} + \hat{M}_{inp}\right)\left(\max_{1 \leq j \leq n}\left[\lambda_i\left(T_i^{(j-1)p} y\right) \langle T_i^{(j-1)p} x - T_i^{jp} x, f'\left(T_i^{(j-1)p} y\right) - f'\left(T_i^{jp} y\right)\rangle\right]\right) \quad (21)$$

with $T_i^{np} x, T_i^{np} y \in A_i$; $\forall y \in A_{i+1}$ for the given $i \in \bar{p}$ with $K(y) < 1$, where $K(y)$ is independent of the particular $T_i^p : \bigcup_{j \in \bar{p}} A_j | A_i \to A_i$ for $i \in \bar{p}$. One gets by using very close arguments to those used in the proof of Theorem 2.4 that $\exists \lim_{n,m \to \infty} D_f\left(T_i^{nmp} x, T_i^{nmp} y\right) = 0$. Then, $\{T_i^{nmp} x\}$ converges to some $z_i \in A_i$ which is proven to be a unique fixed point in the nonempty, closed and convex set $A_i$ for $i \in \bar{p}$. The remaining of the proof is similar to that of Theorem 2.3. The last part of the result follows by applying its first part to each of the $p$ generalized point-dependent $(K, \lambda)$-hybrid $1$-cyclic composite self-mappings $T_i^p : \bigcup_{j \in \bar{p}} A_j | A_i \to A_i$ relative to $D_f$; $\forall i \in \bar{p}$. □

**Remark 2.5**. If $f : D\left|\bigcup_{i \in \bar{p}} A_i \subset X \to (-\infty, \infty]\right.$ is totally convex if it is a continuous strictly convex function which is Gâteaux differentiable in $int\, D$, $dim\, X < \infty$ and $D\left|\bigcup_{i \in \bar{p}} A_i\right.$ is closed, [7]. In view of this result, Theorem 2.3 and Corollary 2.4 are still valid if the condition of its strict total convexity of $f : D\left|\bigcup_{i \in \bar{p}} A_i \to (-\infty, \infty]\right.$ is replaced by its continuity and its strict convexity if the Banach space is finite dimensional. Since $v_f(x, t) \geq \delta_f(t) > 0$; $\forall t \in R_+$, it turns out that if $f : D\left|\bigcup_{i \in \bar{p}} A_i \to (-\infty, \infty]\right.$ is uniformly convex then it is totally convex. Therefore, Theorem 2.3 and Corollary 2.4 still hold if the condition of strict total convexity is replaced with the sufficient one of strict uniform convexity. Note that if a convex function $f$ is totally convex then it is sequentially consistent in the sense that $D_f(x_n, y_n) \to 0$ as $n \to \infty$ if $\|x_n - y_n\| \to 0$ as $n \to \infty$ for any sequences $\{x_n\}$ and $\{y_n\}$ in $D$. □

Some results on weak cluster points of average sequences of relative hybrid self-mappings built from $\{T^n x\}$, for $x \in A_i$ and some $i \in \bar{p}$ are investigated in the following results related to the fixed points of $\{T^n x\}$.

**Theorem 2.6**. Assume that:

1) $X$ is a reflexive space and $f: D \subset X \to (-\infty, \infty]$ be a lower-semicontinuous strictly convex function so that it is Gâteaux differentiable in $int(D)$.

2) A $p$-cyclic self-mapping $T: \bigcup_{i \in \bar{p}} A_i \to \bigcup_{i \in \bar{p}} A_i$ is given defining a composite self-mapping $T^p: \bigcup_{i \in \bar{p}} A_i \to \bigcup_{i \in \bar{p}} A_i$ with $A_i (\neq \emptyset) \subseteq int\, D \subset H$ being bounded and closed; $\forall i \in \bar{p}$, so that its restricted mapping $T_i^p: \bigcup_{j \in \bar{p}} A_j | A_i \to A_i$ to $A_i$, for some given $i \in \bar{p}$, is generalized point-dependent $(1, \lambda_i)$-hybrid relative to $D_f$ for some $\lambda_i: A_i \to \mathbf{R}$ for the given $i \in \bar{p}$ with such a $A_i(\neq \emptyset) \subseteq int\, D \subset H$ being, furthermore, convex.

Define the sequence $\{S_n^{(i)} x\} \equiv \{\frac{1}{n} \sum_{k=0}^{n-1} T^{kp} x\}$ for $x \in A_i$, where $T_i^0 \equiv T^{0p}$ is the identity mapping on $A_i$ so that $T^{0p} x = x$; $\forall x \in A_i$ and assume that $\{T^n x\}$ is bounded for $x \in A_i$. Then, the following properties hold:

**(i)** Every weak cluster point of $\{S_n^{(i)} x\}$ for $x \in A_i$ is a fixed point $v_i \in A_i$ of $T_i^p: \bigcup_{j \in \bar{p}} A_j | A_i \to A_i$ of $T_i^p: \bigcup_{j \in \bar{p}} A_j | A_i \to A_i$ for the given $i \in \bar{p}$.

**(ii)** Define sequences $\{S_n^{(i,j)} x\} \equiv \{\frac{1}{n} \sum_{k=0}^{n-1} T^{kp+j} x\}$ for any integer $1 \leq j \leq p-1$ and $x \in A_i$ where $A_i$ are bounded, closed and convex; $\forall i \in \bar{p}$. Thus, $\{S_n^{(i,j)} x\}$ converges to $v_{i+j} = T^j v_i \in A_{i+j}$ for $x \in A_i$, where $v_i \in A_i$ is a fixed point of $T_i^p: \bigcup_{j \in \bar{p}} A_j | A_i \to A_i$ and a weak cluster point of $\{S_n^{(i)} x\}$ for $x \in A_i$ and $v_{i+j} \in A_{i+j} (1 \leq j \leq p-1)$ is both a fixed point of $T_{i+j}^p: \bigcup_{j \in \bar{p}} A_j | A_{i+j} \to A_{i+j}$ and a weak cluster point of $\{S_n^{(i,j)} x\}$ for $x \in A_i$. Furthermore, $v_{i+j} = T^j v_i$ if $T: \bigcup_{i \in \bar{p}} A_i \to \bigcup_{i \in \bar{p}} A_i$ is continuous. □

## 3. Some extensions with weak convergence to weak cluster points of a class of sequences

Some of the results of Section 2 are now generalized to the case when the subsets of the cyclic mapping do not intersect $T: \bigcup_{i \in \bar{p}} A_i \to \bigcup_{i \in \bar{p}} A_i$, in general, by taking advantage of the fact that best proximity points of such a self-mapping are fixed points of the restricted $T_i^p: \bigcup_{j \in \bar{p}} A_j | A_i \to A_i$ for $i \in \bar{p}$. Weak convergence of averaging sequences to weak cluster points and their links with the best proximity points

in the various subsets of the $p$-cyclic self-mappings is discussed. Firstly, the following result follows from a close proof to that of Theorem 2.6 which is omitted.

**Theorem 3.1**. Let $X$ be a reflexive space and let $f : D \subset X \to (-\infty, \infty]$ be a lower-semicontinuous strictly convex function so that it is Gâteaux differentiable in $int(D)$. Consider the generalized point-dependent $(p \geq 1)$-cyclic hybrid self-mapping $T : \bigcup_{i \in \bar{p}} A_i \to \bigcup_{i \in \bar{p}} A_i$ being $(1, \lambda_i)$ relative to $D_f$ for some $\lambda_i : A_i \to \mathbf{R}$ such that $A_i (\neq \emptyset) \subseteq int\ D \subset H$ are all bounded, convex, closed and with nonempty intersection. Define the sequence $\{S_n x\} \equiv \left\{ \frac{1}{n} \sum_{k=0}^{n-1} T^k x \right\}$ for $x \in \bigcup_{i \in \bar{p}} A_i$, where $T^0$ is the identity mapping on $\bigcup_{i \in \bar{p}} A_i$ and assume that $\{T^n x\}$ is bounded for $x \in \bigcup_{i \in \bar{p}} A_i$. Then, the following properties hold:

**(i)** Every weak cluster point of $\{S_n^{(i)} x\}$ for $x \in A_i$ is a fixed point $v_i \in A_i$ of $T : \bigcup_{i \in \bar{p}} A_i \to \bigcup_{i \in \bar{p}} A_i$.

**(ii)** Define the sequence $\{S_n x\} \equiv \left\{ \frac{1}{n} \sum_{k=0}^{n-1} T^k x \right\}$ for $x \in \bigcup_{i \in \bar{p}} A_i$ which is bounded, closed and convex; $\forall i \in \bar{p}$ and any integer $1 \leq j \leq p-1$. Thus, $\{S_n x\}$ converges to the fixed point $v = Tv \in \bigcap_{i \in \bar{p}} A_i$ of $T : \bigcup_{i \in \bar{p}} A_i \to \bigcup_{i \in \bar{p}} A_i$ for $x \in \bigcup_{i \in \bar{p}} A_i$ which is also a weak cluster point of $\{S_n x\}$. □

**Remark 3.2**. The results of Theorem 2.6 and Theorem 3.1 are extendable without difficulty to the weak cluster points of other related sequences to the considered ones. In particular,

1) Define sequences $\{S_n^{(j)} x\} \equiv \left\{ \frac{1}{n} \sum_{k=0}^{n-1} T^{k+j} x \right\}$, $x \in \bigcup_{i \in \bar{p}} A_i$ for any given finite non-negative integer $j$ under all the hypotheses of Theorem 3.1. With this notation, the sequence considered in such a corollary is $\{S_n x\} \equiv \{S_n^{(0)} x\}$. Direct calculation yields:

$$\left( S_n^{(j)} x - S_n x \right) = \frac{1}{n} \sum_{k=0}^{j-1} \left( T^{k+n} x - T^k x \right) \to 0 \text{ for } x \in \bigcup_{i \in \bar{p}} A_i \text{ as } n \to \infty \text{ since } \left\{ T^{k+n} x - T^k x \right\}, \text{ and then}$$

$\left\{ \sum_{k=0}^{j-1} \left( T^{k+n} x - T^k x \right) \right\}$, is bounded. Then, $S_n^{(j)} x \to v$ weakly which is the same fixed point of $T : \bigcup_{i \in \bar{p}} A_i \to \bigcup_{i \in \bar{p}} A_i$ in $\bigcap_{i \in \bar{p}} A_i$ which is a weak cluster point of $\{S_n^{(j)} x\}$ for $x \in \bigcup_{i \in \bar{p}} A_i$ for any finite non-negative integer $j$.

2) Consider all the hypotheses of Theorem 3.1 and now define sequences $\{S_n^{[j]}x\} \equiv \left\{\frac{1}{n}\sum_{k=0}^{n+j-1}T^k x\right\}$, $x \in \bigcup_{i\in\bar{p}} A_i$ for any given finite non-negative integer $j$. With this notation, the sequence considered in the corollary is $\{S_n x\} \equiv \{S_n^{[0]}x\}$. Direct calculation yields:

$$\left(S_n^{[j]}x - S_n x\right) = \frac{1}{n}\sum_{k=0}^{j-1}\left(T^{k+n}x\right) \to 0$$

weakly for $x \in \bigcup_{i\in\bar{p}} A_i$ as $n \to \infty$ since $\{T^{k+n}x\}$, and then and then $\left\{\sum_{k=0}^{j-1}\left(T^{k+n}x\right)\right\}$, is bounded. Then, $S_n^{[j]}x \to v$ weakly which is the same fixed point of $T: \bigcup_{i\in\bar{p}} A_i \to \bigcup_{i\in\bar{p}} A_i$ in $\bigcap_{i\in\bar{p}} A_i$ which is a weak cluster point of $\{S_n^{[j]}x\}$ for $x \in \bigcup_{i\in\bar{p}} A_i$ and for any finite non-negative integer $j$.

3) Now consider the hypotheses of Theorem 2.6. It turns out that the sequence $\{S_n^{(i,j)}x\} \equiv \left\{\frac{1}{n}\sum_{k=0}^{n-1}T^{kp+j}x\right\}$ for $x \in A_i$ satisfies for any integer $1 \le j \le p-1$

$$S_n^{(i,j)}x = \frac{1}{n+1}\frac{n+1}{n}\left(\sum_{k=0}^{n-1}T^{kp+j}x\right) + \frac{T^j}{n}\left(\sum_{k=0}^{n-1}T^{kp}x - \sum_{k=0}^{n}T^{kp}x\right) = T^j\left(\frac{n+1}{n}S_{n+1}^{(i)}x - \frac{1}{n}T^{np}x\right)$$

$$S_n^{(i,j)}x = \frac{1}{n+1}\frac{n+1}{n}\left(\sum_{k=0}^{n-1}T^{kp}\left(T^j x\right)\right) + \frac{1}{n}\left(\sum_{k=0}^{n-1}T^{kp}\left(T^j x\right) - \sum_{k=0}^{n}T^{kp}x\right)$$

$$= \left(\frac{n+1}{n}S_{n+1}^{(i+j)}x_{i+j} - \frac{1}{n}T^{np}x_{i+j}\right) \to v_{i+j} \equiv T_{i+j}^p v_{i+j} (\in A_{i+k})$$

weakly as $n \to \infty$, where $x_{i+j}(=T^j x) \in A_{i+j}$ since $x \in A_i$, $\{T^n x\}$ is bounded, and

$$\left\{S_{n+1}^{(i,j)}z\right\} \equiv \frac{1}{n}\sum_{k=0}^{n-1}T^{kp}z = \frac{1}{n}\sum_{k=0}^{n-1}\left(T_{i+j}^p\right)^k z$$

for $z \in A_{i+j}$ and $1 \le j \le p-1$. Thus, $v_{i+j}$ is a fixed point of $T_{i+j}^p : \bigcup_{j\in\bar{p}} A_j\big|A_{i+j} \to A_{i+j}$ which is also a weak cluster point of the sequences $\{S_{n+1}^{(i,j)}z\}$ for $1 \le j \le p-1$. However, it is not guaranteed that $v_{i+j} = T^j v_i = T_i^p v_i = T_{i+j}^p v_{i+j}$ without additional hypotheses on $T : \bigcup_{i\in\bar{p}} A_i \to \bigcup_{i\in\bar{p}} A_i$ such as its continuity, or art least, that of the composite mapping $T^j : A_i \to A_{i+j}$ allowing the existence of function of limit equalizing limit of the function at such a fixed point.

4) Now, define $\{S_n^{[i,j]}x\} \equiv \left\{\frac{1}{n}\sum_{k=0}^{n+j-1}T^{kp}x\right\}$ for $x \in A_i$. Note that for $x \in A_i$, $\exists v_i \in A_i$

$$S_n^{[i,j]}x = \frac{1}{n}\sum_{k=0}^{n-1}T^{kp}x + \frac{1}{n}\sum_{k=0}^{j-1}T^{(n+k)p}x = S_n^{(i)}x + \frac{1}{n}\sum_{k=0}^{j-1}T^{(n+k)p}x \to v_i \left(= T_i^p v_i\right)$$

weakly as $n \to \infty$ since $j$ is finite, which is a fixed point in $A_i$ of the composite mapping $T_i^p$ and a weak cluster point of $\left\{S_n^{[i,j]}x\right\}$ for finite $j$. □

By incorporating some background contractive-type conditions for the cyclic self-mapping the above results can be extended to include uniqueness of fixed points as follows.

**Theorem 3.3**. Assume that:

1) Assumption 1 of Theorem 2.6 holds with the restriction of $(X, \|\ \|)$ to be a uniformly convex Banach space.

2) Assumption 2 of Theorem 2.6 holds and, furthermore, all the $p$-cyclic composite mappings with restricted domain $T_i^p : \bigcup_{j \in \overline{p}} A_j | A_i \to A_i$; $\forall i \in \overline{p}$ are either contractive or Meir-Keeler contractions.

Then, the following properties hold

**(i)** Theorem 2.6[(i)-(ii)] holds. Furthermore, each of the mappings $T_i^p : \bigcup_{j \in \overline{p}} A_j | A_i \to A_i$ has a unique fixed point $v_i \in A_i$ which are also best proximity points of $T : \bigcup_{i \in \overline{p}} A_i \to \bigcup_{i \in \overline{p}} A_i$ in $A_i$ so that $v_{i+j} = T^j v_i$; $\forall j \in \overline{p-i}$, $\forall i \in \overline{p}$.

**(ii)** If, in addition, $\bigcap_{i \in \overline{p}} A_i \neq \emptyset$ then, there is a unique fixed point $v \in \bigcap_{i \in \overline{p}} A_i$ of $T : \bigcup_{i \in \overline{p}} A_i \to \bigcup_{i \in \overline{p}} A_i$ and $T_i^p : \bigcup_{j \in \overline{p}} A_j | A_i \to A_i$; $\forall i \in \overline{p}$. □

Theorem 2.9 can be also extended "mutatis-mutandis" to the convergence of weak cluster points of the alternative sequences discussed in Remark 3.2. It is now proven that the sets of fixed points of the restricted composite mapping $T_i^p : \bigcup_{j \in \overline{p}} A_j | A_i \to A_i$; some $i \in \overline{p}$ are convex if such mappings are asymptotically quasi-nonexpansive with respect to $D_f$ in the sense that it has at least a fixed point in $A_i$ and $\limsup_{n \to \infty} D_f\left(v, T_i^{pm} x\right) \leq D_f(v, x)$; $\forall x \in A_i$ and $f : D \subset X \to (-\infty, \infty]$ is a proper strictly convex function. The above concept of asymptotically quasi non-expansive mapping relaxes that of quasi-nonexpansive mapping in [5].

**Theorem 3.4**. Let $f : D \subset X \to (-\infty, \infty]$ be a proper strictly convex function on the Banach space $(X, \|\cdot\|)$ so that is is Gâteaux differentiable in $int\ D$ and consider the restricted composite mapping $T_i^p : \bigcup_{j \in \bar{p}} A_j | A_i \to A_i$ for some $i \in \bar{p}$ built from the $p$-cyclic self-mapping $T : \bigcup_{i \in \bar{p}} A_i \to \bigcup_{i \in \bar{p}} A_i$ so that $A_i$ is nonempty, convex and closed. If $T_i^p : \bigcup_{j \in \bar{p}} A_j | A_i \to A_i$ is quasi-nonexpansive then its set of fixed points $F(T_i^p)$ is a closed and convex subset of $A_i$. □

ACKNOWLEDGEMENTS
The authors thank the Spanish Ministry of Education for its support of this work through Grant DPI2009-07197. They are also grateful to the Basque Government for its support through Grants IT378-10 and SAIOTEK SPE09UN12.